\documentclass[11pt]{article}
\usepackage{mathrsfs}
\usepackage{graphicx}
\usepackage{latexsym,amsmath,amssymb,amsfonts}
\usepackage{color,colordvi,amscd,indentfirst,graphics}
\allowdisplaybreaks

\newtheorem{theorem}{Theorem}[section]

\newtheorem{lemma}[theorem]{Lemma}
\newtheorem{claim}[theorem]{Claim}
\newtheorem{corollary}[theorem]{Corollary}

\def\qed{\hfill \rule{4pt}{7pt}}
\def\pf{\noindent {\it Proof.} }

\begin{document}

\title{Spanning trees with at most $5$ leaves and branch vertices in total of $K_{1,5}$-free graphs}
\author{Pham Hoang Ha\footnote{E-mail address: ha.ph@hnue.edu.vn (Corresponding author).}\\
Department of Mathematics\\
Hanoi National University of Education\\
136 XuanThuy Street, Hanoi, Vietnam\\
\medskip\\
Nguyen Hoang Trang\footnote{E-mail address: trang153@gmail.com.}\\
Foreign Language Speacialized School\\
University of Languages and International Studies\\
Vietnam National University, HaNoi\\
}
\date{}

\maketitle{}

\bigskip

\begin{abstract}
In this paper, we prove that every
$n$-vertex connected $K_{1,5}$-free graph $G$ with $\sigma_4(G)\geq
n-1$ contains a spanning tree with at most $5$ leaves and branch vertices in total. Moreover, the
degree sum condition ``$\sigma_4(G)\geq n-1$" is best possible.
\end{abstract}



\newpage

\section{Introduction}

In this paper, we only consider finite simple graphs. Let $G$ be a
graph with vertex set $V(G)$ and edge set $E(G)$. For any vertex
$v\in V(G)$, we use $N_G(v)$ and $d_G(v)$ (or $N(v)$ and $d(v)$ if
there is no ambiguity) to denote the set of neighbors of $v$ and the
degree of $v$ in $G$, respectively. For any $X\subseteq V(G)$, we
denote by $|X|$ the cardinality of $X$. We define
$N(X)=\bigcup\limits_{x\in X}N(x)$ and $d(X)=\sum\limits_{x\in
X}d(x)$. For an integer $k\geq 1$, we let $N_k(X)=\{x\in V(G)\;|\;|N(x)\cap
X|=k\}$. We
use $G-X$ to denote the graph obtained from $G$ by deleting the
vertices in $X$ together with their incident edges. The subgraph of
$G$ induced by $X$ is denoted by $G[X]$. We define $G-uv$ to be the
graph obtained from $G$ by deleting the edge $uv\in E(G)$, and
$G+uv$ to be the graph obtained from $G$ by adding an edge $uv$
between two non-adjacent vertices $u$ and $v$ of $G$. We write $A:=
B$ to rename $B$ as $A$.

A subset $X\subseteq V(G)$ is called an \emph{independent set} of
$G$ if no two vertices of $X$ are adjacent in $G$.  The maximum size
of an independent set in $G$ is denoted by $\alpha(G)$. For $k\geq
1$, we define $\sigma_k(G)=\min\{\sum\limits_{i=1}^k
d(v_i)\;|\;\{v_1,\ldots,v_k\}$ is an independent set in $G$\}. For
$r\geq 1$, a graph is said to be \emph{$K_{1,r}$-free} if it does
not contain $K_{1,r}$ as an induced subgraph. A $K_{1,3}$-free graph
is also called a \emph{claw-free} graph. 

Let $T$ be a tree. A vertex of degree one is a \emph{leaf} of $T$
and a vertex of degree at least three is a \emph{branch vertex} of
$T$. The set of leaves of $T$
is denoted by $L(T)$ and the set of branch vertices of $T$ is denoted by $B(T)$. For two distinct vertices $u,v$ of $T$, we denote by $P_T[u,v]$
the unique path in $T$ connecting $u$ and $v$ and denote by
$d_T[u,v]$ the distance between $u$ and $v$ in $T$. We define the
\emph{orientation} of $P_T[u,v]$ is from $u$ to $v$. 

There are many known results on the independence
number conditions and the degree sum conditions to ensure that a
connected graph $G$ contains a spanning tree with a bounded number of leaves
or branch vertices.  Win~\cite{Wi79} obtained a
sufficient condition related to the independence number for
$k$-connected graphs having a few leaves, which confirms a conjecture of Las
Vergnas~\cite{LV71}. On the other hand, Broersma and Tuinstra~\cite{BT98} gave a degree
sum condition for a connected graph to contain a spanning tree with
a bounded number of leaves. Beside that, recently, the first named author \cite{H} stated an improvement of Win's result by giving an independence number condition for a graph having a spanning tree which covers a certain subset of $V(G)$ and has at most $l$ leaves. 

In 2012, Kano et al.~\cite{KKMOSY12} presented a degree sum condition for a
connected claw-free graph to have a spanning tree with at most $l$
leaves, which generalizes a result of Matthews and
Sumner~\cite{MS84} and a result of Gargano et al.~\cite{GHHSV04}. Later,
 Chen et al.~\cite{CLX17}, Matsuda et al.~\cite{MOY} and Gould and Shull \cite{GS} also considered the sufficient conditions for a connected
claw-free graph to have a spanning tree with
few leaves or few branch vertices, respectively. 

On the other hand, Kyaw~\cite{Ky09,Ky11} obtained
the sharp sufficient conditions for connected $K_{1,4}$-free graphs to have a spanning tree with few leaves. After that, many researchers also studied sufficient conditions for existence of spanning trees with few leaves or few branch vertices in connected $K_{1,4}$-free graphs (see Chen et al.~\cite{CCH14} and Ha \cite{H1} for examples). 

For the $K_{1,5}$-free graphs, some results were obtained as follows.
\begin{theorem}[{\cite[Chen et al.]{CHH}}]
	Let $G$ be a connected $K_{1,5}$-free graph with $n$ vertices. If
	$\sigma_5(G)\geq n-1$, then $G$ contains a spanning tree with at
	most $4$ leaves.
\end{theorem}
\begin{theorem}[{\cite[Hu and Sun]{HS}}] Let $G$ be a connected $K_{1,5}$-free graph with $n$ vertices.
	If $\sigma_6(G) \geq  n - 1,$ then $G$ contains a spanning tree with at most $5$ leaves.
\end{theorem}

Moreover, many researchers have also studied the degree sum conditions for graphs to have
spanning trees with a bounded number of branch vertices and leaves.
\begin{theorem}[{\cite[Nikoghosyan]{N}, \cite[Saito and Sano]{SS}}]\label{t3}  Let $k\geq 2$ be an integer. If a connected graph $G$ satisfies $\deg_G(x) + \deg_G(y) \geq |G|-k+1$ for every two non-adjacent vertices $x, y \in V (G)$, then $G$ has
	a spanning tree $T$ with $|L(T)| + |B(T)| \leq k + 1.$
\end{theorem}
In 2019, Maezawa et al. improved the previous result by proving the following
theorem.
\begin{theorem}[{\cite[Maezawa et al.]{MMM}}]  Let $k\geq 2$	be an integer. Suppose that a connected graph $G$ satisfies $\max\{\deg_G(x), \deg_G(y)\} \geq \dfrac{|G|-k+1}{2}$ for every two non-adjacent vertices $x, y \in V (G)$, then $G$ has	a spanning tree $T$ with $|L(T)| + |B(T)| \leq k + 1.$
\end{theorem}
Recently, Hanh and the first named author also gave sharp results for the case of claw-free graphs and $K_{1,4}$-free graphs, respectively.
\begin{theorem}[{\cite[Hanh]{Hanh}}] \label{t}
	Suppose that a connected claw-free graph $G$ of order $n$ satisfies $\sigma_{5}(G)\geq n-2.$ Then $G$ has a spanning tree $T$ with $|B(T)|+|L(T)|\leq 5.$ 
\end{theorem}
\begin{theorem}[{\cite[Ha]{H1}}] \label{thm}
	Let $k, m$ be two non-negative intergers {\rm ($ m \leq k+1$)} and let  $G$ be a connected $K_{1,4}$-free graph of order $n$. If $\sigma_{m+2}(G)\geq n-k$, then $G$ has a spanning tree with at most $m+k+2$ leaves and branch vertices. 
\end{theorem}

In this paper, we further consider connected $K_{1,5}$-free graphs.
We give a sufficient condition for a connected $K_{1,5}$-free graph
to have a spanning tree with few leaves and branch vertices in total. More precisely, we prove the following.

\begin{theorem}\label{theo1.9}
Let $G$ be a connected $K_{1,5}$-free graph with $n$ vertices. If
$\sigma_4(G)\geq n-1$, then $G$ contains a spanning tree with at
most $5$ leaves and branch vertices in total.
\end{theorem}

It is easy to see that if a tree has at least $2$ branch vertices then it has at least $4$ leaves. Therefore, we immediately
obtain the following corollary from Theorem~\ref{theo1.9}.

\begin{corollary}\label{coro1.10}
Let $G$ be a connected $K_{1,5}$-free graph with $n$ vertices. If
$\sigma_4(G)\geq n-1$, then $G$ contains a spanning tree with at
most $1$ branch vertices.
\end{corollary}

We end this section by constructing an example to show that the
degree sum condition ``$\sigma_4(G)\geq n-1$" in
Theorems~\ref{theo1.9} is sharp. For an integer $m\geq 1$, let
$D_1,D_2,D_3,D_4$ be vertex-disjoint copies of the
complete graph $K_m$ with $m$ vertices. Let
$xy$ be an edge such that neither $x$ nor $y$ is contained in
$\bigcup\limits_{i=1}^4V(D_i)$. Join $x$ to all the vertices in
$V(D_1)\cup V(D_2)$ and join $y$ to all the vertices in $V(D_3)\cup
V(D_4)$. The resulting graph is denoted by $G$. Then it is easy to
check that $G$ is a connected $K_{1,5}$-free graph with $n=4m+2$
vertices and $\sigma_4(G)=4m=n-2$. However, every spanning tree of
$G$ contains at least $6$ leaves and branch vertices in total.

\section{Proof of the main result}

In this section, we extend the idea of Chen-Ha-Hanh in~\cite{CHH} to prove
Theorem~\ref{theo1.9}.  For this purpose, we need the following
lemma.

\begin{lemma}\label{lem2.1}
Let $G$ be a connected graph such that $G$ does not have a spanning
tree with at most $5$ leaves and branch vertices in total, and let $T$ be a maximal tree of $G$
with $|L(T)|+|B(T)|\in \{6, 7\}$. Then there does not exist a tree $T'$ in $G$ such
that $|L(T')|+|B(T')|\leq 5$ and $V(T')=V(T)$.
\end{lemma}

\pf Suppose for a contradiction that there exists a tree $T'$ in $G$
with at  most $5$ leaves and branch vertices in total and $V(T')=V(T)$. Since $G$ has no spanning
tree with at  most $5$ leaves and branch vertices in total, we see that $V(G)-V(T')\neq\emptyset$.
Hence there must exist two vertices $v$ and $w$ in $G$ such that
$v\in V(T')$ and $w\in N(v)\cap (V(G)-V(T'))$. Let $T_1$ be the tree
obtained from $T'$ by adding the vertex $w$ and the edge $vw$. Then $|L(T_1)|+|B(T_1)|-|L(T')|-|B(T')| \in \{0, 1, 2\}$.

If $|L(T_1)|+|B(T_1)|\in \{6, 7\}$, then $T_1$ contradicts the maximality of
$T$ (since $|V(T_1)|=|V(T)|+1>|V(T)|$). So we may assume that $|L(T_1)|+|B(T_1)|\leq 5$. By repeating this process, we can
recursively construct a set of trees $\{T_i\;|\;i\geq 1\}$ in $G$
such that $|L(T_i)|+|B(T_i)|\leq 5$ and $|V(T_{i+1})|=|V(T_i)|+1$
for each $i\geq 1$. Since $G$ has no spanning tree with at most $5$ leaves and branch vertices in total and $|V(G)|$ is finite, the process must terminate after a
finite number of steps, i.e., there exists some $k\geq 1$ such that
$T_{k+1}$ is a tree in $G$ such that $|L(T_{k+1})|+|B(T_{k+1})|\in \{6, 7\}$. But this contradicts the
maximality of $T$. So the lemma holds.    \qed

\bigskip

\noindent{\bf Proof of Theorem~\ref{theo1.9}.} We prove the theorem
by contradiction. Suppose to the contrary that $G$ contains no
spanning tree with at most $5$ leaves and branch vertices in total. Then every spanning tree of
$G$ contains at least $6$ leaves and branch vertices in total. We choose a maximal tree $T$ of
$G$ with $|L(T)|+|B(T)|\in \{6, 7\}$. 

We consider four cases according to the number of branch vertices
in $T$. (Note that $T$ contains at most two branch vertices.)

\medskip

\textit{Case 1}. $T$ contains two branch vertices and four leaves.

\medskip

Let $s$ and $t$ be the two branch vertices in $T$ and let $U=\{u_1; u_2; u_3; u_4\}$ be the set of leaves of $T.$ Then
$d_T(s)=d_T(t)=3$. Moreover, by the maximality of $T$, we have
$N(U)\subseteq V(T)$. For simplifying notation, let $[k]$ be the set of $\lbrace 1,2,\ldots,k\rbrace$ for some positive integer $k$.

For each $i\in [4]$, let $B_i$ be
the vertex set of the connected component of $T-\{s,t\}$ containing
$u_i$ and let $v_i$ be the unique vertex in $B_i\cap N_T(\{s,t\})$.
Without loss of generality, we may assume that
$\{v_1,v_2\}\subseteq N_T(s)$ and $\{v_3,v_4\}\subseteq N_T(t)$.
For each $1\leq i\leq 4$ and $x\in B_i$, we use $x^-$ and $x^+$ to
denote the predecessor and the successor of $x$ on $P_T[s,u_i]$ or
$P_T[t,u_i]$, respectively (if such a vertex exists). Let $s^+$ be the successor of $s$ on
$P_{T}[s,t]$. Define $P:=V(P_T[s,t])-\{s,t\}$.

For this case, we further choose $T$ such that
\begin{itemize}
\item [$($C1$)$]
$d_T[s,t]$ is as small as possible.
\end{itemize}

\begin{claim}\label{claim2.2}
For all $1\leq i,j\leq 4$ and $i\neq j$, if $x\in N(u_j)\cap B_i$,
then $x\not=u_i,x\not=v_i$ and $x^-\notin N(U-\{u_j\})$.
\end{claim}

\pf Suppose $x=u_i$ or $x=v_i$. Then $T':=T-v_iv_i^-+xu_j$ is a tree
in $G$ with $3$ leaves and $1$ branch vertex such that $V(T')=V(T)$, which contradicts
Lemma~\ref{lem2.1}. So we have $x\not=u_i,x\not=v_i$.

Next, assume $x^-\in N(U-\{u_j\})$. Then there exists some
$k\in [4]-\{j\}$ such that $x^-u_k\in E(G)$. Now,
$T':=T-\{v_iv_i^-,xx^-\}+\{xu_j,x^-u_k\}$ is a tree in $G$ with  $3$ leaves and $1$ branch vertex such that $V(T')=V(T)$, also contradicting Lemma~\ref{lem2.1}. This
proves Claim~\ref{claim2.2}.     \qed

\bigskip

By Claim~\ref{claim2.2}, we know that $U$ is an independent set in
$G$.

\begin{claim}\label{claim2.3}
$N(u_i)\cap P=\emptyset$ for each $i\in [4]$.
\end{claim}

\pf Suppose the assertion of the claim is false. Then there exists
some vertex $x\in P$ such that $xu_i\in E(G)$ for some
$i\in [4]$. Let $T':=T-v_iv_i^{-}+xu_i$, then $T'$ is a tree in $G$
such that $V(T')=V(T)$, $T'$ has $4$ leaves and $2$ branch vertices
$s'$ and $t'$ and
$d_{T'}[s',t']<d_T[s,t]$. But this contradicts the condition (C1). So
the claim holds.   \qed

\begin{claim}\label{claim2.6}
	$N(u_i)\cap \{t\}=\emptyset$ for each $ i\in [2]$.
\end{claim}

\pf Suppose $su_i\in E(G)$ for some $i\in [2]$. Consider the tree $T':=T-v_iv_i^{-}+tu_i$ is a tree in $G$ with
$4$ leaves and $1$ branch vertex such that $V(T')=V(T)$, contradicting Lemma~\ref{lem2.1}.  This proves Claim~\ref{claim2.6}.      \qed

Similarly, we also have 
\begin{claim}\label{claim2.7}
$N(u_i)\cap \{s\}=\emptyset$ for each $3\leq i\leq 4$.
\end{claim}

\begin{claim}\label{claim2.9}
$N_2(U-u_i)\cap B_i=\emptyset$ for each $i\in [4]$. In particular, $N_3(U)=(N_2(U)-N(u_i))\cap B_i=\emptyset$ for each $i\in [4].$
\end{claim}
\pf For the sake of convenience, we may assume by symmetry that
$i\in [2]$.

Suppose this is false. Then there exists some vertex $x\in
(N_2(U-u_i))\cap B_i$ for some $i\in [2]$. By applying
Claim~\ref{claim2.2}, we have $x\not=u_i$ and $ x\not=v_i.$ 

Since $x\in
N_2(U-u_i)\cap B_i$ there
must exist two distinct indices $j,k\in [4]-\{i\}, j < k,$ such that
$xu_j,xu_k\in E(G)$. Set

$$T':=\left\{\begin{array}{ll}T-\{v_jv_j^{-},v_kv_k^{-}\}+\{xu_j,xu_k\},
& \;\mbox{ if } j=3-i,\\
T-\{ss^{^+},v_kv_k^{-}\}+\{xu_j,xu_k\},
& \;\mbox{ if } 3\leq j< k \leq 4,
\end{array}\right.$$
Then
$T'$ is a tree in $G$
with $1$ branch vertex and $4$ leaves such that $V(T')=V(T)$, contradicting Lemma~\ref{lem2.1}.
 \qed
\bigskip

By Claims~\ref{claim2.2} and \ref{claim2.9},  $\{u_i\}$, $N(u_i)\cap B_i$, and 
$(N(U-\{u_i\})\cap B_i)^-$ are
pairwise disjoint subsets in $B_i$ for each $i\in [4]$ (where
$(N(U-\{u_i\})\cap B_i)^-=\{x^-\;|\;x\in N(U-\{u_i\})\cap B_i\}$) and $N_3(U)=(N_2(U)-N(u_i))\cap B_i=\emptyset$ for
each $i\in [4]$.
Then for each $i\in [4]$, we conclude that
\begin{align*}
|B_i|&\geq 1+|N(u_i)\cap B_i|+|(N(U-\{u_i\})\cap
B_i)^-|\\
&= 1+|N(u_i)\cap B_i|+|N(U-\{u_i\})\cap B_i|\\
&=1+\sum\limits_{j=1}^4|N(u_j)\cap B_i|.
\end{align*}
By applying
Claim~\ref{claim2.3}, we obtain
\begin{align*}
\sum\limits_{i=1}^4|N(u_i)\cap P|=0.
\end{align*}
On the other hand, by Claims \ref{claim2.6}-\ref{claim2.7} we obtain that
\begin{align*}\sum\limits_{i=1}^{4}|N(u_i)\cap \{s\}|\leq 2, \sum\limits_{i=1}^{4}|N(u_i)\cap \{t\}|\leq 2.\end{align*}
Note that $N(U)\subseteq V(T)$. Now, we conclude
that
\begin{align*}
|V(T)|&=\sum\limits_{i=1}^{4}|B_i|+|V(P_{T}[s,t])|\\
&\geq\sum\limits_{i=1}^{4}\left(\sum\limits_{j=1}^4|N(u_j)\cap
B_i|+1\right)+
\left(\sum\limits_{i=1}^{4}|N(u_i)\cap \{s\}|+ \sum\limits_{i=1}^{4}|N(u_i)\cap \{t\}|-2+\sum\limits_{i=1}^{4}|N(u_i)\cap P|\right)\\
&=2+\sum\limits_{i=1}^{4}\sum\limits_{j=1}^{4}|N(u_j)\cap
B_i|+\sum\limits_{i=1}^{4}|N(u_i)\cap \{s,t\}|+\sum\limits_{i=1}^{4}|N(u_i)\cap P|\\
&=\sum\limits_{j=1}^{4}|N(u_j)\cap V(T)|+2\\
&=\sum\limits_{j=1}^{4}d(u_j)+2\\
&=d(U)+2.
\end{align*}
Since $U$ is an independent set in $G$, we have
\begin{align*}
n-1\leq \sigma_4(G)\leq d(U)\leq |V(T)|-2\leq n-2,
\end{align*}
a contradiction.

\medskip
\textit{Case 2}. $T$ contains two branch vertices and five leaves.

\medskip

Let $s$ and $t$ be the two branch vertices in $T$ such that
$d_T(s)=4$ and $d_T(t)=3$. Let $U=\{u_1; u_2; u_3; u_4; u_5\}$ be the set of leaves of $T.$ For each $i\in [5]$, let $B_i$ be
the vertex set of the connected component of $T-\{s,t\}$ containing
$u_i$ and let $v_i$ be the unique vertex in $B_i\cap N_T(\{s,t\})$.
Without loss of generality, we may assume that
$\{v_1,v_2,v_3\}\subseteq N_T(s)$ and $\{v_4,v_5\}\subseteq N_T(t)$.
For each $i\in [5]$ and $x\in B_i$, we use $x^-$ and $x^+$ to
denote the predecessor and the successor of $x$ on $P_T[s,u_i]$ or
$P_T[t,u_i]$, respectively (if such a vertex exists). Let $s^+$ and
$t^-$ be the successor of $s$ and the predecessor of $t$ on
$P_{T}[s,t]$, respectively. Define $P:=V(P_T[s,t])-\{s,t\}$.

For this case, we choose $T$ such that
\begin{itemize}
	\item [$($D1$)$]
	$d_T[s,t]$ is as small as possible, and
	\item [$($D2$)$]
	 $\sum\limits_{i=1}^3 |B_i|$ is as large as
	possible, subject to $($D1$)$.
\end{itemize}

\begin{claim}\label{claim3.2}
	For all $1\leq i,j\leq 5$ and $i\neq j$, if $x\in N(u_j)\cap B_i$,
	then $x\not=u_i,x\not=v_i$ and $x^-\notin N(U-\{u_j\})$.
\end{claim}

\pf Suppose $x=u_i$ or $x=v_i$. Then $T':=T-v_iv_i^-+xu_j$ is a tree
in $G$ with $4$ leaves and at most $2$ branch vertices such that $V(T')=V(T)$. Then this contradicts
either Lemma~\ref{lem2.1} or the proof of Case
1. So we have $x\not=u_i,x\not=v_i$.

Next, assume $x^-\in N(U-\{u_j\})$. Then there exists some
$k\in [5]-\{j\}$ such that $x^-u_k\in E(G)$. Now,
$T':=T-\{v_iv_i^-,xx^-\}+\{xu_j,x^-u_k\}$ is a tree in $G$ with $4$
leaves and at most $2$ branch vertices such that $V(T')=V(T)$, also contradicting either Lemma~\ref{lem2.1} or the proof of Case
1. This
proves Claim~\ref{claim3.2}.     \qed

\bigskip

By Claim~\ref{claim3.2}, we know that $U$ is an independent set in
$G$. Since $G$ is $K_{1,5}$-free, we have $N_5(U)=\emptyset$.

\begin{claim}\label{claim3.3}
	$N(u_i)\cap P=\emptyset$ for each $4\leq i\leq 5$.
\end{claim}

\pf Suppose the assertion of the claim is false. Then there exists
some vertex $x\in P$ such that $xu_i\in E(G)$ for some
$i\in\{4,5\}$. Let $T':=T-tv_i+xu_i$, then $T'$ is a tree in $G$
with $5$ leaves such that $V(T')=V(T)$, $T'$ has two branch vertices
$s$ and $x$, $d_{T'}(s)=4$, $d_{T'}(x)=3$ and
$d_{T'}[s,x]<d_T[s,t]$. But this contradicts the condition (D1). So
the claim holds.   \qed

\begin{claim}\label{claim3.4}
	If $P\neq \emptyset$, then $\sum\limits_{i=1}^3 |N(u_i)\cap
	\{x\}|\leq 1$ for each $x\in P$.
\end{claim}

\pf Suppose to the contrary that there exists some vertex $x\in P$
such that $\sum\limits_{i=1}^3 |N(u_i)\cap \{x\}|\geq 2$. Then there
exist two distinct indices $j,k\in [3]$ such that $xu_j,xu_k\in E(G)$.
Let $T':=T-\{sv_j,sv_k\}+\{xu_j,xu_k\}$, then $T'$ is a tree in $G$
with $5$ leaves such that $V(T')=V(T)$, $T'$ has two branch vertices
$x$ and $t$, $d_{T'}(x)=4$, $d_{T'}(t)=3$ and
$d_{T'}[x,t]<d_T[s,t]$, contradicting the condition (D1). This
completes the proof of Claim~\ref{claim3.4}.     \qed

\begin{claim}\label{claim3.6}
	$N(u_i)\cap \{s\}=\emptyset$ for each $4\leq i\leq 5$.
\end{claim}

\pf Suppose $su_i\in E(G)$ for some $i\in\{4,5\}$. If $P=\emptyset$,
then we have $st\in E(T)$ and $T':=T-st+su_i$ is a tree in $G$ with
$|L(T')|+|B(T')|=5$ and $V(T')=V(T)$, contradicting Lemma~\ref{lem2.1}. So we
may assume that $P\neq\emptyset$ and hence $s^+\neq t$. By applying
Claims~\ref{claim3.2} and~\ref{claim3.3}, we deduce that
$N(u_i)\cap\{s^+,v_1,v_2,v_3\}=\emptyset$.

Suppose that $s^+v_j\in E(G)$ for some $j\in [3]$. Then
$T':=T-\{ss^+,sv_j\}+\{su_i,s^+v_j\}$ is a tree in $G$ with $4$
leaves and $2$ branch vertices such that $V(T')=V(T).$ By the same argument as in the proof of Case
1, we can derive a contradiction. So we
conclude that $N(s^+)\cap \{v_1,v_2,v_3\}=\emptyset$.

Now, assume there exits two distinct $j,k\in [3]$ such that
$v_jv_k\in E(G)$. Then by Claim~\ref{claim3.2}, we see that $u_k\neq
v_k$. Let $T':=T-\{sv_j,tv_i\}+\{su_i,v_jv_k\}$, then $T'$ is a tree
in $G$ with $2$ branch vertices and $5$ leaves such that $V(T')=V(T)$, $T'$ has two branch
vertices $s$ and $v_k$, $d_{T'}(s)=4$, $d_{T'}(v_k)=3$ and
$d_{T'}[s,v_k]<d_T[s,t]$, contradicting the condition (D1).
Therefore, $v_1,v_2$ and $v_3$ are pairwise non-adjacent in $G$.

But then, $\{s^+,u_i,v_1,v_2,v_3\}$ is an independent set and
$G[\{s,s^+,u_i,v_1,v_2,v_3\}]$ is an induced $K_{1,5}$ of $G$, again
a contradiction. This proves Claim~\ref{claim3.6}.      \qed

\begin{claim}\label{claim3.7}
$\sum\limits_{i=1}^5 |N(u_i)\cap \{t\}|\leq 3$.
\end{claim}

\pf Suppose for a contradiction that $\sum\limits_{i=1}^5 |N(u_i)\cap \{t\}|\geq 4.$\\

 If $P=\emptyset$ then we have
$st\in E(G)$. Since $\sum\limits_{i=1}^5 |N(u_i)\cap \{t\}|\geq 4$,
there exists some $j\in [3]$ such that $tu_j\in E(G)$. Let
$T':=T-st+tu_j$, then $T'$ is a tree in $G$ with $4$ leaves and $2$ branch vertices such that $V(T')=V(T).$ Repeating the same argument as in the proof of Case
1, we can deduce a contradiction. 

Otherwise, $P\not=\emptyset$, then $t^{-}\not=s.$ It follows
from Claim~\ref{claim3.3} that $N(u_i)\cap\{t^-\}=\emptyset$ for
each $4\leq i\leq 5$. Suppose that $t^-u_j\in E(G)$ for some
$j\in [3]$. Since $t\in N_4(U)$, there exists some
$k\in [3]-\{j\}$ such that $tu_k\in E(G)$. Let
$T':=T-\{sv_j,tt^-\}+\{tu_k,t^-u_j\}$, then $T'$ is a tree in $G$
with $4$ leaves and $2$ branch vertices such that $V(T')=V(T)$, which contradicts
the proof of Case 1. Therefore, we deduce that $N(U)\cap
\{t^-\}=\emptyset$. But then, $(N(t)\cap U)\cup\{t^-\}$ is an
independent set and $G[(N(t)\cap U)\cup\{t,t^-\}]$ is an induced
$K_{1,5}$ of $G$, again a contradiction. So the claim
holds.    \qed

\begin{claim}\label{claim3.8}
We have $N_3(U-\{u_i\})\cap B_i=\emptyset$ for every $i\in [5].$ In particular, we obtain	$N_4(U)=\emptyset$.
\end{claim}

\pf Suppose to the contrary that there exists some vertex $x\in
N_3(U-\{u_i\})$ for some $i\in [5]$. By
Claim~\ref{claim3.2}, we know that $x^-, x^+\notin N_3(U-\{u_i\})$. 

Suppose that $x^-x^+\in E(G)$. Since $x\in N_3(U-\{u_i\})$, there
must exist two distinct $j,k\in [5]-\{i\}$ such that
$xu_j,xu_k\in E(G)$. Then
$T':=T-\{v_jv_j^-,xx^-,xx^+\}+\{xu_j,xu_k,x^-x^+\}$ is a tree in $G$
with $4$ leaves and at least $2$ branch vertices such that $V(T')=V(T)$, contradicting either Lemma~\ref{lem2.1} or the proof of Case 1.
Hence $x^-x^+\notin E(G)$.

Then $(N(x)\cap
(U-\{u_i\}))\cup\{x^-,x^+\}$ is an independent set and $G[(N(x)\cap
(U-\{u_i\})\cup\{x,x^-, x^+\}]$ is an induced $K_{1,5}$ of $G$, contradicting the
assumption that $G$ is $K_{1,5}$-free. This completes the proof of
Claim~\ref{claim3.8}. \qed

\begin{claim}\label{claim3.9}
We have	$(N_3(U)-N(u_i))\cap B_i=\emptyset$ for each $1\leq i\leq 5$.
\end{claim}

\pf Suppose this is false. Then there exists some vertex $x\in
(N_3(U)-N(u_i))\cap B_i$ for some $1\leq i\leq 5$. By applying
Claim~\ref{claim3.2}, we have $x\not=u_i, x\not=v_i$ and
$x^-,x^+\notin N(U-\{u_i\})$.

Suppose that $x^-x^+\in E(G)$. Since $x\in N_3(U)-N(u_i)$, there
must exist two distinct indices $j,k\in [5]-\{i\}$ such that
$xu_j,xu_k\in E(G)$. Then
$T':=T-\{v_jv_j^-,xx^-,xx^+\}+\{xu_j,xu_k,x^-x^+\}$ is a tree in $G$
with $4$ leaves and at least $2$ branch vertices such that $V(T')=V(T)$, contradicting either Lemma~\ref{lem2.1} or the proof of Case 1.
Hence $x^-x^+\notin E(G)$.

Now, $(N(x)\cap U)\cup\{x^-,x^+\}$ is an independent set and
$G[(N(x)\cap U)\cup\{x,x^-,x^+\}]$ is an induced $K_{1,5}$ of $G$,
giving a contradiction. So the assertion of the claim holds.   \qed

\begin{claim}\label{claim3.10}
	$N(u_j)\cap B_i=\emptyset$ for all $4\leq i\leq 5$ and $1\leq j\leq
	3$. In particular, $N_3(U)\cap N(u_i)\cap B_i=\emptyset$ for each
	$4\leq i\leq 5$.
\end{claim}

\pf Suppose the assertion of the claim is false. Then there exists
some vertex $x\in B_i$ such that $xu_j\in E(G)$ for some
$i\in\{4,5\}$ and $j\in [3]$. By Claim~\ref{claim3.2}, we have
$x\not=u_i$ and $ x\not=v_i$. Let $T':=T-xx^-+xu_j$, and let $B_k'$ be the
vertex set of the connected component of $T'-\{s,t\}$ containing
$u_k$ for each $1\leq k\leq 3$. It is easy to check that $T'$ is a
tree in $G$ with $5$ leaves such that $V(T')=V(T)$, $T'$ has two
branch vertices $s$ and $t$, $d_{T'}(s)=4$, $d_{T'}(t)=3$,
$d_{T'}[s,t]=d_T[s,t]$ and $\sum\limits_{k=1}^3
|B_k'|=\sum\limits_{k=1}^3 |B_k|+|V(P_T[x,u_i])|>\sum\limits_{k=1}^3
|B_k|$. But this contradicts the condition (D2). This proves
Claim~\ref{claim3.10}. \qed

\begin{claim}\label{claim3.11}
	$|N_3(U)\cap N(u_i)\cap B_i|\leq 1$ for each $1\leq i\leq 3$.
\end{claim}

\pf Suppose for a contradiction that there exist two distinct
vertices $x,y\in N_3(U)\cap N(u_i)\cap B_i$ for some
$i\in [3]$. Without loss of generality, we may assume that
$x\in V(P_T[s,y])$. By Claim~\ref{claim3.2}, we have
$x\not=u_i, x\not=v_i$,$y\not=u_i, y\not=v_i$, $x^-\notin N(U)$ and $x^+\notin
N(U-\{u_i\})$. In particular, $x^+\neq y$. Since $x,y\in N_3(U)\cap
N(u_i)$, there exist two distinct $j,k\in [5]-\{i\}$ such
that $xu_j,yu_k\in E(G)$. We may assume that $x^-x^+,x^+u_i\notin
E(G)$; for otherwise,
$$T':=\left\{\begin{array}{ll}T-\{sv_i,xx^-,xx^+,yy^+\}+\{xu_i,xu_j,x^-x^+,yu_k\},
& \;\mbox{ if } x^-x^+\in E(G),\\
T-\{sv_i,xx^+,yy^-\}+\{xu_j,x^+u_i,yu_k\}, & \;\mbox{ if } x^+u_i\in
E(G),
\end{array}\right.$$
is a tree in $G$ with $4$ leaves and $2$ branch vertices such that $V(T')=V(T).$ By the same argument as in the proof of Case
1, we can deduce a contradiction. But then, $(N(x)\cap U)\cup\{x^-,x^+\}$ is an
independent set and $G[(N(x)\cap U)\cup\{x,x^-,x^+\}]$ is an induced
$K_{1,5}$ of $G$, again a contradiction. So the claim holds.
\qed

\begin{claim}\label{claim3.12}
	For each $1\leq i\leq 3$, if $u_iv_i\in E(G)$, then $N_3(U)\cap
	N(u_i)\cap B_i=\emptyset$.
\end{claim}

\pf Suppose to the contrary that $u_iv_i\in E(G)$ and there exists
some vertex $x\in N_3(U)\cap N(u_i)\cap B_i$ for some
$i\in [3]$. By Claim~\ref{claim3.2}, we have $x\neq v_i$. Since
$x\in N_3(U)\cap N(u_i)$, there exists some
$j\in[5]-\{i\}$ such that $xu_j\in E(G)$. Let
$T':=T-\{sv_i,xx^-\}+\{u_iv_i,xu_j\}$, then $T'$ is a tree in $G$
with $4$ leaves and two branch vertices such that $V(T')=V(T).$ By the same argument as in the proof of Case
1, we give a contradiction. This completes the proof of Claim~\ref{claim3.12}.     \qed

\begin{claim}\label{claim3.13}
	For each $1\leq i\leq 3$, if $su_i\in E(G)$, then $N_3(U)\cap
	N(u_i)\cap B_i=\emptyset$.
\end{claim}

\pf For the sake of convenience, we may assume by symmetry that
$i=1$. Suppose the assertion of the claim is false. Then there
exists some vertex $x\in N_3(U)\cap N(u_1)\cap B_1$. By applying
Claims~\ref{claim3.2} and~\ref{claim3.12}, we know that
$x\notin\{u_1,v_1\}$ and $N(u_1)\cap \{v_1,v_2,v_3\}=\emptyset$.

Suppose $v_1v_j\in E(G)$ for some $j\in\{2,3\}$. Then
$T':=T-\{sv_1,sv_j\}+\{su_1,v_1v_j\}$ is a tree in $G$ with $4$
leaves and $2$ branch vertices such that $V(T')=V(T).$ By the same argument as in the proof of Case
1, we can deduce a contradiction. So we
have $v_1v_2,v_1v_3\notin E(G)$.

Next, assume that $v_2v_3\in E(G)$. Then $u_2\neq v_2$ and $u_3\neq
v_3$ by Claim~\ref{claim3.2}. If there exists some $j\in\{2,3\}$
such that $xu_j\in E(G)$, then $T':=T-\{sv_2,sv_3\}+\{v_2v_3,xu_j\}$
is a tree in $G$ with $4$ leaves and $2$ branch vertices such that $V(T')=V(T)$, contradicting Case 1. Hence $xu_2,xu_3\notin E(G)$. Then, since $x\in
N_3(U)\cap N(u_1)$, we conclude that $xu_4,xu_5\in E(G)$. Let
$T':=T-\{sv_2,tt^-,xx^-\}+\{su_1,v_2v_3,xu_4\}$. If $P=\emptyset$,
then $t^-=s$, and $T'$ is a tree in $G$ with $4$ leaves and $2$ branch vertices such that
$V(T')=V(T)$, giving a contradiction with the proof of Case 1. So we deduce that
$P\neq\emptyset$. But then, $T'$ is a tree in $G$ with $5$ leaves
such that $V(T')=V(T)$, $T'$ has two branch vertices $s$ and $v_3$,
$d_{T'}(s)=4$, $d_{T'}(v_3)=3$ and $d_{T'}[s,v_3]=1<d_T[s,t]$,
contradicting the condition (D1). Therefore, $v_1,v_2$ and $v_3$ are
pairwise non-adjacent in $G$.

We now consider the vertex $s^+$. We will show that
$N(s^+)\cap\{u_1,v_1,v_2,v_3\}=\emptyset$.

We first prove that $s^+u_1\notin E(G)$. Suppose this is false. Let $T':=T-\{ss^+\}+\{s^+u_1\}$, then $T'$ is a tree in $G$ with $4$
leaves and $2$ branch vertices such that $V(T')=V(T).$ By the same argument as in the proof of Case
1, we can deduce a contradiction.

Finally, we show that $s^+v_2,s^+v_3\notin E(G)$. Suppose not, and
let $s^+v_j\in E(G)$ for some $j\in\{2,3\}$. If there exists some
$k\in\{4,5\}$ such that $xu_k\in E(G)$, then
$T':=T-\{ss^+,sv_j\}+\{s^+v_j,xu_k\}$ is a tree in $G$ with $4$
leaves and $2$ branch vertices such that $V(T')=V(T).$ Repeating the same argument as in the proof of Case
1, we can deduce a contradiction. Therefore, we have $xu_4,xu_5\notin E(G)$. Since $x\in N_3(U)\cap
N(u_1)$, we deduce that $xu_2,xu_3\in E(G)$. Let
$T':=T-\{ss^+,sv_j,xx^-,xx^+\}+\{su_1,s^+v_j,xu_2,xu_3\}$, then $T'$
is a tree in $G$ with $4$ leaves and $2$ branch vertices such that $V(T')=V(T)$, again a
contradiction. Hence $N(s^+)\cap\{u_1,v_1,v_2,v_3\}=\emptyset$.

Now, $\{s^+,u_1,v_1,v_2,v_3\}$ is an independent set and
$G[\{s,s^+,u_1,v_1,v_2,v_3\}]$ is an induced $K_{1,5}$ of $G$,
giving a contradiction. So the assertion of the claim holds. \qed

\bigskip

By Claim~\ref{claim3.2}, $\{u_i\}$, $N(u_i)\cap B_i$,
$(N(U-\{u_i\})\cap B_i)^-$ and $(N_2(U)-N(u_i))\cap B_i$ are
pairwise disjoint subsets in $B_i$ for each $i\in [5]$, where
$(N(U-\{u_i\})\cap B_i)^-=\{x^-\;|\;x\in N(U-\{u_i\})\cap B_i\}$.
Recall that $N_5(U)=N_4(U)=(N_3(U)-N(u_i))\cap B_i=\emptyset$ (for
each $1\leq i\leq 5$) by Claims~\ref{claim3.8} and~\ref{claim3.9}.
Then for each $i\in [3]$, we conclude that
\begin{align*}
|B_i|&\geq 1+|N(u_i)\cap B_i|+|(N(U-\{u_i\})\cap
B_i)^-|+|(N_2(U)-N(u_i))\cap B_i|\\
&= 1+|N(u_i)\cap B_i|+|N(U-\{u_i\})\cap B_i|+|(N_2(U)-N(u_i))\cap
B_i|\\
&=1+\sum\limits_{j=1}^5|N(u_j)\cap B_i|-|N_3(U)\cap N(u_i)\cap
B_i|\\
&\geq\sum\limits_{j=1}^5|N(u_j)\cap B_i|+|N(u_i)\cap\{s\}|,
\tag{$1$}
\end{align*}
where the last inequality follows from Claims~\ref{claim3.11}
and~\ref{claim3.13}. Similarly, for each $4\leq i\leq 5$, we have
\begin{align*}
|B_i|&\geq 1+|N(u_i)\cap B_i|+|(N(U-\{u_i\})\cap
B_i)^-|+|(N_2(U)-N(u_i))\cap B_i|\\
&= 1+|N(u_i)\cap B_i|+|N(U-\{u_i\})\cap B_i|+|(N_2(U)-N(u_i))\cap
B_i|\\
&=1+\sum\limits_{j=1}^5|N(u_j)\cap B_i|-|N_3(U)\cap N(u_i)\cap
B_i|\\
&=1+\sum\limits_{j=1}^5|N(u_j)\cap B_i|+|N(u_i)\cap\{s\}|, \tag{$2$}
\end{align*}
where the last equality follows from Claims~\ref{claim3.6}
and~\ref{claim3.10}.

For each $1\leq i\leq 5$, we define $d_i=|N(u_i)\cap P|$. Then
$d_4=d_5=0$ by Claim~\ref{claim3.3}. By applying
Claim~\ref{claim3.4}, we know that $N(u_1)\cap P, N(u_2)\cap P$ and
$N(u_3)\cap P$ are pairwise disjoint. Therefore,
\begin{align*}
|P|\geq \sum\limits_{i=1}^5 d_i=\sum\limits_{i=1}^5|N(u_i)\cap P|.
\end{align*}

By combining Claim~\ref{claim3.7}, we have 
\begin{align*}
|V(P_T[s,t])|=2+|P|\geq \sum\limits_{i=1}^5|N(u_i)\cap\{t\}|+\sum\limits_{i=1}^5|N(u_i)\cap
P| - 1. \tag{$3$}
\end{align*}

Note that $N(U)\subseteq V(T)$. By (1), (2) and (3), we conclude
that
\begin{align*}
|V(T)|&=\sum\limits_{i=1}^{3}|B_i|+\sum\limits_{i=4}^{5}|B_i|+|V(P_{T}[s,t])|\\
&\geq\sum\limits_{i=1}^{3}\left(\sum\limits_{j=1}^5|N(u_j)\cap
B_i|+|N(u_i)\cap\{s\}|\right)+
\sum\limits_{i=4}^{5}\left(1+\sum\limits_{j=1}^5|N(u_j)\cap
B_i|+|N(u_i)\cap\{s\}|\right)\\
&\ \ \ \
+\left(\sum\limits_{i=1}^5|N(u_i)\cap\{t\}|+\sum\limits_{i=1}^5|N(u_i)\cap
P|-1\right)\\
&=2+\sum\limits_{i=1}^{5}\sum\limits_{j=1}^{5}|N(u_j)\cap
B_i|+\sum\limits_{i=1}^{5}|N(u_i)\cap \{s,t\}|+\sum\limits_{i=1}^{5}|N(u_i)\cap P|\\
&=\sum\limits_{j=1}^{5}|N(u_j)\cap V(T)|+1\\
&=\sum\limits_{j=1}^{5}d(u_j)+1\\
&=d(U)+1.
\end{align*}
Since $U$ is an independent set in $G$, we have
\begin{align*}
n-1\leq \sigma_4(G)\leq \sigma_5(G)-1 \leq d(U)-1\leq |V(T)|-2\leq n-2,
\end{align*}
a contradiction.

\textit{Case 3}. $T$ contains one branch vertex and five leaves.

\medskip

Let $r$ be the unique branch vertex in $T$ with $d_T(r)=5$ and let
$N_T(r)=\{v_1,v_2,v_3,v_4,v_5\}$. For each $i\in [5]$, let $B_i$ be
the vertex set of the connected component of $T-\{r\}$ containing
$u_i$ and let $v_i$ be the unique vertex in $B_i\cap N_T(\{r\})$. For each $i\in [5]$ and $x\in B_i$, we use $x^-$ and $x^+$ to
denote the predecessor and the successor of $x$ on $P_T[r,u_i]$, respectively (if such a vertex exists).\\

Since $G$ is $K_{1,5}$-free, there
exist two distinct indices $i,j\in [5]$ such that $v_iv_j\in E(G)$.
Let $T':=T-rv_i+v_iv_j$. If $v_j$ is a leaf of $T$, then $T'$ is a
tree in $G$ with $4$ leaves and $1$ branch vertex such that $V(T')=V(T)$, which contradicts
Lemma~\ref{lem2.1}. So we may assume that $v_j$ has degree two in
$T$. Then $T'$ is a tree in $G$ with $5$ leaves such that
$V(T')=V(T)$, $T'$ has two branch vertices $r$ and $v_j$,
$d_{T'}(r)=4$ and $d_{T'}(v_j)=3$. By the same argument as in the
proof of Case 2, we can also derive a contradiction.

\medskip

\textit{Case 4}. $T$ contains one branch vertex and six leaves.

\medskip

Let $r$ be the unique branch vertex in $T$ with $d_T(r)=6$ and let $U=\{u_i\}_{i=1}^{6}$ be the set of leaves of $T.$ For each $i\in [6]$, let $B_i$ be
the vertex set of the connected component of $T-\{r\}$ containing
$u_i$ and let $v_i$ be the unique vertex in $B_i\cap N_T(\{r\})$. For each $i\in [6]$ and $x\in B_i$, we use $x^-$ and $x^+$ to
denote the predecessor and the successor of $x$ on $P_T[r,u_i]$, respectively (if such a vertex exists).

\begin{claim}\label{claim4.14}
For all $1\leq i,j\leq 6$ and $i\neq j$, if $x\in N(u_j)\cap B_i$,
then $x\not=u_i, x\not=v_i$ and $x^-\notin N(U-\{u_j\})$.
\end{claim}

\pf Suppose $x=u_i$ or $x=v_i$. Then $T':=T-v_iv_i^-+xu_j$ is a tree
in $G$ with $5$ leaves and $1$ branch vertex such that $V(T')=V(T)$. By the same argument as in the proof of Case
3, we can derive a contradiction. So we have $x\not=u_i,x\not=v_i$.

Next, assume $x^-\in N(U-\{u_j\})$. Then there exists some
$k\in [6]-\{j\}$ such that $x^-u_k\in E(G)$. Now,
$T':=T-\{v_iv_i^-,xx^-\}+\{xu_j,x^-u_k\}$ is a tree in $G$ with $5$
leaves and $1$ branch vertex such that $V(T')=V(T).$ By the same argument as in the proof of Case
3, we can deduce a contradiction. This
proves Claim~\ref{claim4.14}.     \qed

By applying Claim~\ref{claim4.14}, we deduce that $U$ is an
independent set in $G$.

\begin{claim}\label{claim4.15}
For every $1\leq i, j \leq 6$ and $i\neq j,$ if $v_iv_j \in E(G)$ then $N(u_i)\cap B_k=\emptyset$ and $N(u_k)\cap B_i=\emptyset$ for each $k \in [6]- \{i, j\}$.
\end{claim}

\pf Suppose to the contrary that there exists some vertex $x\in B_k$
such that $xu_i\in E(G)$. By Claim ~\ref{claim4.14} we have $x\not= u_k.$ Let
$T':=T-\{v_iv_i^-, v_jv_j^-\}+ \{xu_i, v_iv_j\}$. Then $T'$ is a tree in $G$
with $5$ leaves such that $V(T')=V(T)$, $T'$ has two branch
vertices $x,r.$ This implies a contradiction by using the proof of Case 2. 

Now, suppose that there exists some vertex $x\in B_i$
such that $xu_k\in E(G)$. By Claim ~\ref{claim4.14} we have $x\not= u_i.$ Let
$T':=T-\{v_iv_i^-, v_jv_j^-\}+ \{xu_k, v_iv_j\}$. Then $T'$ is a tree in $G$
with $5$ leaves such that $V(T')=V(T)$, $T'$ has two branch
vertices $x,r.$ This also gives a contradiction by using the same arguments as in the proof of Case 2. 

The proof of Claim \ref{claim4.15} is completed.
 \qed
 
Since $G$ is $K_{1,5}$-free, there
exist two distinct indices $i,j\in [5]$ such that $v_iv_j\in E(G)$. Without loss of generality, we may assume that $v_1v_2 \in E(G).$ 

Set $U_1=\{u_1, u_2, u_3, u_4\}.$ By Claim \ref{claim4.15} we obtain that 
$$N(U_1)\cap B_j = N(\{u_i\}_{i=1}^2)\cap B_j\,\, \text{for all } j \in \{1, 2\},$$ 
and
$$N(U_1)\cap B_j = N(\{u_i\}_{i=3}^4)\cap B_j\,\, \text{for all } j \in \{3, 4, 5, 6\}.$$

By Claim~\ref{claim4.14}, $\{u_i\}$, $N(u_i)\cap B_i$, and
$(N(U_1-\{u_i\})\cap B_i)^-$ are
pairwise disjoint subsets in $B_i$ for each $i\in [4]$, where
$(N(U_1-\{u_i\})\cap B_i)^-=\{x^-\;|\;x\in N(U_1-\{u_i\})\cap B_i\}$.
Recall that $N_4(U_1)=(N_3(U_1)-N(u_i))\cap B_i= (N_2(U_1)-N(u_i))\cap B_i= \emptyset$ (for
each $1\leq i\leq 4$). 
Then for each $i\in [4]$, we conclude that
\begin{align*}
|B_i|&\geq 1+|N(u_i)\cap B_i|+|(N(U_1-\{u_i\})\cap
B_i)^-|\\
&= 1+|N(u_i)\cap B_i|+|N(U_1-\{u_i\})\cap B_i|+|(N_2(U_1)-N(u_i))\cap
B_i|\\
&=1+\sum\limits_{j=1}^4|N(u_j)\cap B_i|.  \tag{$6$}
\end{align*}

On the other hand, by Claim~\ref{claim4.14}, $\{u_i\}$, $N(u_3)\cap B_i$, and
$(N(u_4)\cap B_i)^-$ are
pairwise disjoint subsets in $B_i$ for each $i\in \{5, 6\}$.
Then for each $i\in \{5, 6\}$, we conclude that
\begin{align*}
|B_i|&\geq 1+|N(u_3)\cap B_i|+|(N(u_4)\cap
B_i)^-|\\
&= 1+|N(u_3)\cap B_i|+|N(u_4)\cap B_i|\\
&=1+\sum\limits_{j=1}^4|N(u_j)\cap B_i|.  \tag{$7$}
\end{align*}
By (6) and (7), we conclude
that
\begin{align*}
|V(T)|&=1+\sum\limits_{i=1}^{6}|B_i|\\
&\geq 1+\sum\limits_{i=1}^{6}\left(1+\sum\limits_{j=1}^6|N(u_j)\cap
B_i|\right)\\
&=7+\sum\limits_{i=1}^{6}\sum\limits_{j=1}^{4}|N(u_j)\cap
B_i|\\
&\geq 3+\sum\limits_{i=1}^{6}\sum\limits_{j=1}^{4}|N(u_j)\cap
B_i|+\sum\limits_{j=1}^{4}|N(u_j)\cap \{r\}|\\
&=\sum\limits_{j=1}^{4}|N(u_j)\cap V(T)|+3\\
&=\sum\limits_{j=1}^{4}d(u_j)+3\\
&=d(U_1)+3.
\end{align*}
Since $U_1$ is an independent set in $G$, we have
\begin{align*}
n-1\leq \sigma_4(G)\leq  d(U_1)\leq |V(T)|-3\leq n-3.
\end{align*}
This also gives a contradiction.

This completes the proof of Theorem~\ref{theo1.9}. \qed

\end{document}